\documentclass[a4paper,11pt]{amsart}
\usepackage{a4wide}
\usepackage{amsmath,amsthm}
\usepackage[psamsfonts]{amssymb}
\usepackage[utf8]{inputenc}
\usepackage{enumerate}
\usepackage{xcolor,graphicx}
\usepackage[hyphens]{url}
\usepackage[ddmmyyyy,hhmmss]{datetime}
\usepackage{hyperref}

\allowdisplaybreaks[4]                  

\newcommand{\field}[1]{\mathbb{{#1}}}
\newcommand{\ideal}[1]{\mathfrak{{#1}}}

\newcommand{\R}{\field{R}}

\newcommand{\Q}{\field{{Q}}}
\newcommand{\K}{\field{{K}}}
\renewcommand{\L}{\field{{L}}}
\newcommand{\p}{\ideal{p}}
\renewcommand{\a}{\ideal{a}}
\newcommand{\Norm}{\textrm{\upshape N}}
\newcommand{\Cl}{\mathcal C\!\ell}
\newcommand{\dd}{\,\mathrm{d}}
\newcommand{\W}{\mathcal W}

\DeclareMathOperator{\Imm}{Im}
\DeclareMathOperator{\Ree}{Re}
\DeclareMathOperator{\Ch}{Ch}
\DeclareMathOperator{\Sh}{Sh}
\DeclareMathOperator{\odd}{odd}
\DeclareMathOperator{\supp}{supp}
\DeclareMathOperator{\Gal}{Gal}
\DeclareMathOperator{\dilog}{dilog}

\newcounter{fixedfig}
\newenvironment{fixedfig}{
\refstepcounter{fixedfig}
\centerline\bgroup%
\def\caption##1{\textsc{Figure \arabic{fixedfig}}: ##1}%
\begin{tabular}{c}%
}
{%
\end{tabular}%
\egroup%
}

\newtheorem{corollary}{Corollary}
\newtheorem{proposition}{Proposition}
\newtheorem{theorem}{Theorem}

\theoremstyle{remark}
\newtheorem*{acknowledgments}{Acknowledgments}

\begin{document}
\title[Breaking the 4 barrier]
{Breaking the 4 barrier for the bound of a generating set of the class group}

\author[L.~Greni\'{e}]{Lo\"{\i}c Greni\'{e}}
\address[L.~Greni\'{e}]{Dipartimento di Ingegneria gestionale, dell'informazione e della produzione\\
         Universit\`{a} di Bergamo\\
         viale Marconi 5\\
         I-24044 Dal\-mi\-ne\\
         Italy}
\email{loic.grenie@gmail.com}

\author[G.~Molteni]{Giuseppe Molteni}
\address[G.~Molteni]{Dipartimento di Matematica\\
         Universit\`{a} di Milano\\
         via Saldini 50\\
         I-20133 Milano\\
         Italy}
\email{giuseppe.molteni1@unimi.it}

\keywords{Class group, generators of the class group}
\subjclass[2020]{Primary 11R04, 11R29; Secondary 11Y40}

\date{\today.
}

\begin{abstract}
Let $\K$ be a field of degree $n$ and discriminant with absolute value $\Delta$. Under the assumption
of the validity of the Generalized Riemann Hypothesis, we provide a new algorithm to compute a set of
generators of the class group of $\K$ and prove that the norm of the ideals in that set is $\leq
(4-1/(2n))\log^2\Delta$, except for a finite number of fields of degree $n\leq 4$. For those fields,
the conclusion holds with the slightly larger limit $(4-1/(2n)+1/(2n^2))\log^2\Delta$. When the
cardinality of $\Cl$ is odd the bounds improve to $(4-2/(3n))\log^2\Delta$, again with finitely many
exceptions in degree $n\leq 4$, and to $(4-2/(3n)+3/(8n^2))\log^2\Delta$ without exceptions.
\end{abstract}

\maketitle

\begin{center}
To appear in Math. Comp. 2025
\end{center}

\section{Introduction}\label{sec:1A}
Let $\K/\Q$ be a number field, with $n=[\K:\Q]$ the degree and $r_1$, $r_2$ the number of real and pair of
complex embeddings of $\K$ (so that $n=r_1+2r_2$). Let $\Delta$ be the absolute value of the absolute
discriminant of $\K$. Let $\Cl$ be the class group of $\K$, i.e. the quotient of the group of
fractional ideals in $\K$ by the subgroup of principal fractional ideals.\\
It is known that $\Cl$ is a finite and abelian group. Buchmann's algorithm is an efficient method to
compute $\Cl$, but it needs as basic ingredient a list of generators. Let $T_\K$ be the minimum of
integers $T$ such that $\{[\p]\colon \p\text{ prime, } \Norm\p \leq T\}$ is a generating set for
$\Cl$.\\
The classical result of Minkowski shows that $T_\K \leq c(r_1,r_2)\sqrt{\Delta}$ for a suitable constant
$c(r_1,r_2)$ depending on $r_1$ and $r_2$ in an explicit way. The value of $c(r_1,r_2)$ has been improved
several time, with significative contributions by Rogers~\cite{Rogers}, Mulholland~\cite{Mulholland} and
Zimmert~\cite{Zimmert} in general, and by de la Maza~\cite{Maza} for certain signatures in degree $n\leq
10$.

In spite of these improvements, the dependence of these bounds on the square root of the discriminant
greatly reduces the possibility to use them in computations. The situation improves dramatically under
the assumption of the validity of the Generalized Riemann Hypothesis. In fact, Bach~\cite{Bach2} proved
that in this case $T_\K \leq c\log^2\Delta$ with $c=12$ in all cases, and $c=4+o(1)$ when $\Delta$
diverges. Recently the authors have proved~\cite{GrenieMolteni7} that already $c = 4.01$ suffices and that
one can take $T_\K \leq 4(\log\Delta + \log\log\Delta)^2$ (other and more precise upper bounds are offered
in that paper).\\
The constant $4$ appearing in these bounds resisted -- so far -- to any attempt to reduce it.
Nevertheless, Belabas, Diaz y Diaz and Friedman~\cite{BelabasDiazyDiazFriedman} modified Bach's approach
and obtained an algorithm producing an upper bound of $T_\K$ for every fixed $\K$; the algorithm has been
improved in~\cite{GrenieMolteni7}. In all tests the latter algorithm has been able to produce a bound
which is essentially $T_\K \leq \log^2\Delta$: this suggests that the constant $4$ represents only a
barrier for the current technology of proof, but that probably the true value of this parameter is
considerably smaller; as small as $1$, may be, or even smaller: see~Figure~\ref{fig:2A}.

In this paper we show that in fact it is possible to reduce the constant below $4$.
\begin{theorem}\label{th:1A}
Assume GRH. Let $\K$ be a field with degree $n$ and absolute value of the discriminant $\Delta$. If
either $n\geq 5$ or $n\leq 4$ and $\Delta\leq\Delta_{-}$ or $\Delta\geq\Delta_{+}$ with $\Delta_{\pm}$ as
in Table~\ref{tab:1A}, one has:
\begin{equation}\label{eq:1A}
T_\K \leq \Big(4-\frac{1/2}{n}\Big) \log^2\Delta.
\end{equation}
\end{theorem}
\begin{table}[h]
\def\z{\phantom{0}}
\def\m#1{\multicolumn{1}{c}{#1}}
\begin{tabular}{crrr}
 $n$              & \m{2} & \m{3} & \m{4} \\
  \hline
 $\log\Delta_{-}$ &    63 &   353 &   833 \\
 $\log\Delta_{+}$ & 12184 & 18455 & 27911
\end{tabular}
\vspace{0.2cm}
\caption{Values of $\Delta_{\pm}$ for Theorem~\ref{th:1A}.}\label{tab:1A}
\end{table}
The value of $\Delta_-$ is already quite large and probably adequate to cover most numerical
applications. Nevertheless, the procedure can be used to deduce explicit values for $\Delta_\alpha$ such
that the bound $T_\K \leq (4-\alpha)\log^2\Delta$ holds whenever $\Delta \geq \Delta_\alpha$, for every
$\alpha\in[0,1/(2n)]$ and for some $\alpha$ all exceptions disappear. As an example we prove the
following result.
\begin{theorem}\label{th:2A}
Assume GRH. Let $\K$ be a field with degree $n$ and absolute value of the discriminant $\Delta$. Then
\[
T_\K \leq \Big(4 - \frac{1/2}{n} + \frac{1/2}{n^2}\Big) \log^2\Delta.
\]
\end{theorem}
\smallskip

We can also adapt the argument for Theorem~\ref{th:1A} to show that with a lower number of ideals it is
still possible to generate at least a large subgroup of $\Cl$. In fact, let $\Cl^2$ be the subgroup of
squares of elements of $\Cl$. It coincides with $\Cl$ if and only if $\Cl$ has odd order and contains the
subgroup $\Cl_{\odd}$ of elements in $\Cl$ having odd order in any case.
Identifying $\Cl$ with $\Gal(\L/\K)$, the Galois group of the Hilbert class field $\L$ for $\K$, the
subgroup $\Cl^2$ corresponds (via Galois correspondence) to the compositum of all quadratic and unramified
extensions of $\K$: this fact allows to appreciate that $\Cl^2$ is usually a large part of $\Cl$.

Let $T'_\K$ be the minimum of integers $T$ such that $\{[\p]\colon \p\text{ prime, } \Norm\p \leq T\}$ is
a generating set for $\Cl^2$.
\begin{theorem}\label{th:3A}
Assume GRH. Let $\K$ be a field with degree $n$ and absolute value of the discriminant $\Delta$. If
either $n\geq 5$ or $n\leq 4$ and $\Delta\leq\Delta'_{-}$ or $\Delta\geq\Delta'_{+}$ with
$\Delta'_{\pm}$ as in Table~\ref{tab:2A}, then
\begin{equation}\label{eq:2A}
T'_\K \leq \Big(4-\frac{2/3}{n}\Big) \log^2\Delta.
\end{equation}
\end{theorem}
\begin{table}[h]
\def\z{\phantom{0}}
\def\m#1{\multicolumn{1}{c}{#1}}
\begin{tabular}{crrrr}
 $n$               & \m{2} & \m{3} & \m{4} \\
  \hline
 $\log\Delta'_{-}$ &    65 &   379 &   993 \\
 $\log\Delta'_{+}$ &  7968 &  9832 & 13108
\end{tabular}
\vspace{0.2cm}
\caption{Values of $\Delta'_{\pm}$ for Theorem~\ref{th:3A}.}\label{tab:2A}
\end{table}
\noindent %
And similarly to Theorem~\eqref{th:2A} we also have a bound for $T'_\K$ which holds true without exceptions.
\begin{theorem}\label{th:4A}
Assume GRH. Let $\K$ be a field with degree $n$ and absolute value of the discriminant $\Delta$. Then
\[
T'_\K \leq \Big(4 - \frac{2/3}{n} + \frac{3/8}{n^2}\Big) \log^2\Delta.
\]
\end{theorem}
Hence, if the cardinality of $\Cl$ is odd, Theorems~\ref{th:3A}--\ref{th:4A} improve directly the
conclusion in Theorems~\ref{th:1A}--\ref{th:2A}. This remark is not trivial since there are cases where
the parity of the class number is known in advance, without the need of a full computation of $\Cl$: for
example, according to a celebrated result of Weber this is what happens for $\Q[2^k]^+$, the maximal real
subfield of the cyclotomic field of order $2^k$, for every $k$. Other cases have been classified:
see~\cite{Berger}, \cite{Hasse4} and~\cite{ConnerHurrelbrink}.
\medskip

Actually, the reduction stated in Theorems~\ref{th:1A}--\ref{th:4A} is quite small and becomes smaller
and smaller as the degree of $\K$ increases, but it does not degrade with the discriminant.

The general strategy for the proofs is similar to the one used in~\cite{GrenieMolteni7}, the new result
is made possible by a new way to deal with the support of the test function. In particular, the level $L$
of its support is not assumed to coincide with (the logarithm of) the level $T$ as usually done, but it
is allowed to be larger. This introduces the need to deal with a certain sum on ideals, but it also
allows more flexibility. A judicious choice of $L$ in terms of $T$ produces the conclusion for
Theorems~\ref{th:1A}--\ref{th:2A}. For Theorems~\ref{th:3A}--\ref{th:4A} we further take advantage of the
special structure of $\Cl^2$.
\smallskip

The new strategy can be adapted to also produce a new algorithm for the computation of $T_\K$. The new
algorithm reduces the size of the bound generally by a factor around 13\% with respect to the algorithm
which is currently used in PARI/GP~\cite{PARI2} for the same purpose. Section~\ref{sec:6A} contains a
presentation of this new algorithm and an analysis of its performance: under this point of view,
Theorems~\ref{th:1A}--\ref{th:2A} double as a quality assessment for the new algorithm.

\medskip

\begin{acknowledgments}
Experiments presented in this paper were carried out using the PlaFRIM experimental testbed, supported by
Inria, CNRS (LABRI and IMB), Universit\'{e} de Bordeaux, Bordeaux INP and Conseil R\'{e}gional
d’Aquitaine~\cite{BORDEAUX-CLUSTER}. %
The authors are members of the INdAM group GNSAGA.
\end{acknowledgments}

\section{Initial computations}\label{sec:2A}
Let $\W$ be the set of even functions $F\colon \R\to\R$ such that:
\begin{enumerate}[a)]
\item $F$ is continuous;
\item $\exists\varepsilon>0$ such that the function $F(x)e^{(\frac{1}{2}+\varepsilon)x}$ is integrable
      and of bounded variation;
\item $(F(0)-F(x))/x$ is of bounded variation.
\end{enumerate}
For any $F\in\W$ let $\phi(s):=\!\int_\R F(x)e^{(s-1/2)x}\dd x$. Let $\chi$ be any character for $\Cl$.
Then, Weil--Poitou's explicit formula~\cite{Poitou} states that
\begin{align}
\delta_{\chi}(\phi(1) + \phi(0)) - \sum_{\rho_\chi} \phi(\rho_\chi)
 =& \sum_{\a} \frac{\Lambda_{\K}(\a)}{\sqrt{\Norm\a}}(\chi(\a)+\overline{\chi(\a)}) F(\log \Norm\a)    \notag\\
& + F(0)((\gamma+\log 8\pi)n - \log \Delta)
  - I(F) n
  + J(F) r_1,                                                                                          \label{eq:3A}
\end{align}
with
\[
I(F):= \int_0^{+\infty} \frac{F(0) {-} F(x)}{2\Sh(x/2)}\dd x
\quad,\quad
J(F):= \int_0^{+\infty} \frac{F(x)}{2\Ch(x/2)}\dd x
\]
and where $\delta_{\chi}$ is $1$ if $\chi$ is the trivial character $\chi_0$, $0$ otherwise, and
$\rho_\chi$ describes the set of critical zeros of $L(s,\chi)$.
Note that all terms appearing on the second line do not depend on the character.

Let $T < T_\K$. Then the prime ideals with norm $\leq T$ generate a proper subgroup of $\Cl$. Hence there
exists a character $\chi\neq \chi_0$ of $\Cl$ which is trivial on the subgroup. Suppose that $F$ is
supported in $[-L,L]$ with $L \geq \log T$. Subtracting the formulas for $L(s,\chi_0)=\zeta_\K$ and
$L(s,\chi)$ produces the equality:
\[
\phi(1) + \phi(0) - \sum_{\rho_\K} \phi(\rho_\K) + \sum_{\rho_\chi} \phi(\rho_\chi)
= \sum_{\a} \frac{\Lambda_{\K}(\a)}{\sqrt{\Norm\a}}|1-\chi(\a)|^2 F(\log \Norm\a),
\]
since $2-\chi(\a)-\overline{\chi(\a)}=|1-\chi(\a)|^2$. Ideals $\a$ with $\Norm\a \leq T$ do not
contribute to the sum, since for them $\chi(\a)=1$. Moreover, the von Mangoldt function selects ideals
which are of the kind $\p^m$, where $\p$ is a prime ideal. Assuming $e^{L/2} < T$, i.e. $L<2\log T$, the
condition $\Norm\p^m \leq e^L$ forces $\Norm\p < T$ whenever $m\geq 2$, so that $\chi(\p^m)=\chi(\p)^m=1$
also in this case. Therefore, under the assumption that $\log T\leq L < 2\log T$ the formula actually
says that
\begin{equation}\label{eq:4A}
\phi(1) + \phi(0) - \sum_{\rho_\K} \phi(\rho_\K) + \sum_{\rho_\chi} \phi(\rho_\chi)
= \sum_{\substack{T< \Norm\p \leq e^L}} \frac{\Lambda_{\K}(\p)}{\sqrt{\Norm\p}}|1-\chi(\p)|^2 F(\log \Norm\p),
\end{equation}
where the sum now runs on prime ideals. Suppose that $\phi(1/2+it)\geq 0$ for every $t\in\R$ and assume
GRH: then $\sum_{\rho_\chi} \phi(\rho_\chi) \geq 0$, and, since $|1-\chi(\p)|\leq 2$, the formula implies
that
\begin{equation}\label{eq:5A}
4\int_0^{+\infty}F(x)\Ch(x/2)\dd x
 =    \phi(1) + \phi(0)
 \leq \sum_{\rho_\K} \phi(\rho_\K)
    + 4\sum_{\substack{T< \Norm\p \leq e^L}} \frac{\Lambda_{\K}(\p)}{\sqrt{\Norm\p}} F(\log \Norm\p).
\end{equation}
Therefore, in case
\begin{equation}\label{eq:6A}
4\int_0^{+\infty}F(x)\Ch(x/2)\dd x
> \sum_{\rho_\K} \phi(\rho_\K)
    + 4\sum_{\substack{T< \Norm\p \leq e^L}} \frac{\Lambda_{\K}(\p)}{\sqrt{\Norm\p}} F(\log \Norm\p)
\quad \text{then} \quad T_\K \leq T.
\end{equation}
Since $\phi(1/2+it) = \int_\R F(x)e^{itx}\dd x = \hat{F}(t)$, we can meet the condition $\phi(1/2+it)\geq
0$ setting $F= \psi\ast\psi$ where $\psi$ is real, even, stepwise $C^1$ and supported in $[-L/2,L/2]$.
In this case one has
\begin{align*}
\int_0^{\infty}F(x)\Ch(x/2)\dd x
=2 \Big[\int_{0}^{+\infty}\psi(x)\Ch(x/2)\dd x\Big]^2,
\end{align*}
and $\hat{F} = \hat{\psi}^2$. Thus, in terms of $\psi$ the criterion in~\eqref{eq:6A} becomes:
\begin{equation}\label{eq:7A}
8 \Big[\int_{0}^{+\infty} \psi(x)\Ch(x/2)\dd x\Big]^2
   > \sum_{\gamma_\K} \big|\hat{\psi}(\gamma_\K)\big|^2
    + 4\sum_{\substack{T< \Norm\p \leq e^L}} \frac{\Lambda_{\K}(\p)}{\sqrt{\Norm\p}} \psi\ast\psi(\log \Norm\p)
\implies
T_\K \leq T,
\end{equation}
where $\gamma_\K :=\Imm \rho_\K$. We further specialize the function setting $\psi = \psi^+ +
\psi^-$, where $\psi^-(x):=\psi^+(-x)$ and $\psi^+$ is positive, stepwise $C^1$ and supported in
$[0,L/2]$.
Then $\hat{\psi}(t) = 2\Ree[\widehat{\psi^+}(t)]$ so that $|\hat{\psi}(t)|^2 \leq 4
|\widehat{\psi^+}(t)|^2 = 4\widehat{\psi^+\ast\psi^-}(t)$, and the criterion in~\eqref{eq:7A} becomes:
\begin{multline}\label{eq:8A}
2\Big[\int_{0}^{L/2} \psi^+(x)\Ch(x/2) \dd x \Big]^2\\
> \sum_{\gamma_\K} \widehat{\psi^{+}{\ast}\psi^{-}}(\gamma_\K)
    + \sum_{\substack{T< \Norm\p \leq e^L}} \frac{\Lambda_{\K}(\p)}{\sqrt{\Norm\p}} \psi\ast\psi(\log \Norm\p)
\implies
T_\K \leq T.
\end{multline}
The key point here is that the series on zeros in~\eqref{eq:8A} involves values of a Fourier trasform, so
that we can apply one more time the explicit formula~\eqref{eq:3A} for $\zeta_\K$ with $F =
\psi^+\ast\psi^-$ to write $\sum_{\gamma_\K} \widehat{\psi^{+}{\ast}\psi^{-}}(\gamma_\K)$ as a sum on
ideals, getting that:
if there exists $\psi^+$ which is positive, stepwise $C^1$, supported in $[0,L/2]$ with
\begin{align*}
2\Big[\int_{0}^{L/2} &\psi^+(x)\Ch(x/2) \dd x \Big]^2
>  4 \int_0^{+\infty} \psi^+\ast\psi^-(x)\Ch(x/2)\dd x                                                         \\
&
   + (\psi^+\ast\psi^-)(0)(\log \Delta - (\gamma+\log 8\pi)n)
   + I(\psi^+\ast\psi^-) n
   - J(\psi^+\ast\psi^-) r_1                                                                                   \\
&
   - 2\sum_{\a} \frac{\Lambda_{\K}(\a)}{\sqrt{\Norm\a}} \psi^+\ast\psi^-(\log\Norm\a)
   +  \sum_{\substack{T< \Norm\p \leq e^L}} \frac{\Lambda_{\K}(\p)}{\sqrt{\Norm\p}} \psi\ast\psi(\log \Norm\p),
\end{align*}
then $T_\K \leq T$.
Since
\[
2\int_0^{+\infty} \psi^+\ast\psi^-(x) \Ch(x/2) \dd x
 = \int_{0}^{L/2} \psi^+(w)e^{w/2} \dd w \int_{0}^{L/2} \psi^+(w) e^{-w/2} \dd w,
\]
the test can be written also as
\begin{align}
2\Big[\int_{0}^{L/2} &\psi^+(x)\Sh(x/2) \dd x \Big]^2                                                     \notag \\
>&  (\psi^+\ast\psi^-)(0)(\log \Delta - (\gamma+\log 8\pi)n)
   + I(\psi^+\ast\psi^-) n
   - J(\psi^+\ast\psi^-) r_1                                                                              \notag \\
&  - 2\sum_{\a} \frac{\Lambda_{\K}(\a)}{\sqrt{\Norm\a}} \psi^+\ast\psi^-(\log\Norm\a)
   +  \sum_{\substack{T< \Norm\p \leq e^L}} \frac{\Lambda_{\K}(\p)}{\sqrt{\Norm\p}} \psi\ast\psi(\log \Norm\p)
\quad
\implies T_\K \leq T.                                                                                     \label{eq:9A}
\end{align}
The assumptions ensure that $\psi^+\ast\psi^-(x) \geq 0$ for every $x$ so that the sum $-2\sum_\a \dots$
is negative but very small in absolute value, because the support of $\psi^+\ast\psi^-$ is $[-L/2,L/2]$
so that the sum ranges only up to $e^{L/2}$. We retain this term in the next computations but in the
final proofs it will be estimated with $0$ but for the cases $n \leq 5$. Moreover, $\psi\ast\psi =
\psi^+\ast\psi^+ + 2\psi^+\ast\psi^- + \psi^-\ast\psi^- $. Functions $\psi^+\ast\psi^-$ and
$\psi^-\ast\psi^-$ are supported in $[-L/2,L/2]$ and $[-L,0]$, respectively, therefore they do not
contribute to the second sum ranging in $(T,e^L]$ (recall that we are assuming $T> e^{L/2}$).

Now we have to make a choice for $\psi^+$. For this purpose we notice that the inequality
in~\eqref{eq:9A} essentially reads
\[
1
   > \frac{\int_0^{L/2} ({\psi^+}(x))^2\dd x}{2\big(\int_{0}^{L/2} \psi^+(x)\Sh(x/2) \dd x\big)^2}(\log\Delta - (\gamma+\log(8\pi))n)
        + \text{ lower order terms}.
\]
To optimize the conclusion we look for a function $\psi^+$ minimizing the coefficient of $\log\Delta$.
Cauchy--Schwarz inequality shows that
\[
\frac{\int_0^{L/2} ({\psi^+}(x))^2\dd x}{2\big(\int_{0}^{L/2} \psi^+(x)\Sh(x/2) \dd x\big)^2}
\geq \frac{1}{\Sh(L/2)-L/2},
\]
with equality only if $\psi^+(x)$ is proportional to $\Sh(x/2)$. Thus, a possible choice would be to set
$\psi^+(x)$ equal to $e^{x/2} - e^{-x/2}$ in $[0,L/2]$ and $0$ otherwise. We have tested this choice, and
we have noticed that the term $e^{-x/2}$ does not significantly change the coefficient of $\log\Delta$,
but it increases secondary terms that prevent the possibility to deduce the second part of
Theorem~\ref{th:1A} and Theorem~\ref{th:2A}. Hence, also in order to simplify the next computations, we
take $\psi^+(x) = e^{x/2}$ supported in $[0,L/2]$. This $\psi^+$ yields:

\begin{align*}
&\psi^-(x) = e^{-x/2}  \text{ on $[-L/2,0]$},   \\
&\psi(x)   = e^{|x|/2} \text{ on $[-L/2,L/2]$}, \\
&2\int_0^{+\infty} \psi^+(x) \Sh(x/2)\dd x = e^{L/2} - 1 - L/2
\intertext{and}
&\psi^+*\psi^-(x)
 = \int_{[x-L/2,x]\cap[-L/2,0]} \!\!\!\!\!\!\!\!\!\!\!\!\!\!\!\!\!\!\!\!\!\!\!\!\!\! e^{(x-u)/2} e^{-u/2}\dd u
 = \begin{cases}
    0                      & \text{if $x<-L/2$}      ,\\
    e^{(L+x)/2}-e^{-x/2}   & \text{if $x\in[-L/2,0]$},\\
    e^{(L-x)/2}-e^{ x/2}   & \text{if $x\in[0,L/2]$} ,\\
    0                      & \text{if $x>L/2$}.
   \end{cases}
\end{align*}
Moreover,
%
%
\begin{align*}
2I(\psi^+*\psi^-)
%
&= e^{L/2} L - 2(e^{L/2}-1)\log(e^{L/2}-1)
   + 4(e^{L/2}-1)\log 2
\intertext{and}
2J(\psi^+*\psi^-)
&= 2\int_0^{L/2} \frac{e^{L/2}-e^{x}}{e^{x}+1}\dd x        \\
&=  e^{L/2} L - 2(e^{L/2}+1)\log(e^{L/2}+1) + 2(e^{L/2}+1)\log 2.
\end{align*}
Finally,
\begin{align*}
\psi^+\ast\psi^+(x)
&= \int_{[x-L/2,x]\cap[0,L/2]} \!\!\!\!\!\!\!\!\!\!\!\!\!\!\!\!\!\!\! e^{(x-u)/2}e^{u/2} \dd u
 = \begin{cases}
    0            & \text{if $x<0$}        ,\\
    e^{x/2}x     & \text{if $x\in[0,L/2]$},\\
    e^{x/2}(L-x) & \text{if $x\in[L/2,L]$},\\
    0            & \text{if $x>L$}.
   \end{cases}
\end{align*}
For future use we also notice that
\begin{equation}\label{eq:10A}
F(x)
 = \begin{cases}
        e^{|x|/2}(|x|-2) + 2e^{(L-|x|)/2} & \text{if $|x|\in[0,L/2]$},\\
        e^{|x|/2}(L-|x|)                  & \text{if $|x|\in[L/2,L]$},\\
        0                                 & \text{if $|x| > L$},
       \end{cases}
\end{equation}
and
%
\begin{subequations}
\begin{align}
I(F)
  &= 2\big(e^{L/2}-1\big)\log\Big(\frac{4}{1-e^{-L/2}}\Big) - \frac{L^2}{4} + L
   - \frac{\pi^2}{6} - 2\dilog(-e^{-L/2}),                                      \label{eq:11aA}\\
J(F)
  &= 2\big(e^{L/2}+1\big)\log\Big(\frac{2}{1+e^{-L/2}}\Big) + \frac{L^2}{4} - L
   - \frac{\pi^2}{12} - 2\dilog(-e^{-L/2}) + \dilog(-e^{-L}).                   \label{eq:11bA}
\end{align}
\end{subequations}
%
%
With the formulas above, the criterion~\eqref{eq:9A} becomes
\begin{align}
\big[e^{L/2} &- 1 - L/2\big]^2
>  2(e^{L/2}-1)\log \Delta                                                             \notag \\
&  - (e^{L/2} L - 2(e^{L/2}+1)\log(e^{L/2}+1) + 2(e^{L/2}+1)\log 2) r_1                \notag \\
&  - (- e^{L/2} L + 2(e^{L/2}-1)\log(e^{L/2}-1) + 2(e^{L/2}-1)(\gamma+\log 2\pi)) n    \notag \\
&  + 4\sum_{\Norm\a\leq e^{L/2}} \Lambda_{\K}(\a)\Big(1 - \frac{e^{L/2}}{\Norm\a}\Big)
   + 2\sum_{\substack{T< \Norm\p \leq e^L}} \log(\Norm\p)\log(e^L/\Norm\p)
\quad
\implies T_\K \leq T.                                                                  \label{eq:12A}
\end{align}
Note that only the formula for $\psi^+\ast\psi^+(x)$ in $[L/2,L]$ matters here, since we have assumed
$T > e^{L/2}$. \\
The coefficients of $r_1$ and $n$ are essentially linear in $e^{L/2}$, hence we introduce a couple of
functions $\alpha$ and $\beta$ via the identities
\begin{align*}
 \alpha(e^L)e^{L/2} &:= \phantom{-} e^{L/2} L - 2(e^{L/2}+1)\log(e^{L/2}+1) + 2(e^{L/2}+1)\log 2            ,\\
  \beta(e^L)e^{L/2} &:=           - e^{L/2} L + 2(e^{L/2}-1)\log(e^{L/2}-1) + 2(e^{L/2}-1)(\gamma+\log 2\pi).
\end{align*}
Functions $\alpha$ and $\beta$ are increasing, positive and bounded.
%
%
We further simplify a bit the inequality noticing that
\[
\big[e^{L/2} - 1 - L/2\big]^2 \geq e^L - Le^{L/2} - 2 e^{L/2}.
\]
Introducing this lower bound in~\eqref{eq:12A} we get a weaker (but simpler) version of the test, saying
that
\begin{align}
e^L
>\,&  e^{L/2}(2\log \Delta + 2 - 2e^{-L/2}\log \Delta - \alpha(e^L)\,r_1 - \beta(e^L)\,n +L)    \notag \\
   &\!+ 4 \sum_{\Norm\a\leq e^{L/2}} \Lambda_{\K}(\a)\Big(1 - \frac{e^{L/2}}{\Norm\a}\Big)
      + 2\sum_{\substack{T< \Norm\p \leq e^L}} \Lambda_{\K}(\p)\log(e^L/\Norm\p)
\quad
\implies T_\K \leq T.                                                                           \label{eq:13A}
\end{align}
Now we have to choose $L$ in terms of $T$. The second sum on ideals is asymptotically estimated by a
constant times $(e^L-T)\log(e^L/T)$, by the Prime Ideal Theorem. Hence, using this estimation, the
inequality in the test essentially says:
\[
e^L
>  2e^{L/2}\log \Delta
   + \Omega\big((e^L-T)\log(e^L/T)\big)
   + \text{ lower order terms}.
\]
In particular we need $e^L \gg (e^L-T)\log(e^L/T)$ which forces $e^L/T \ll 1$. Thus, we set $e^L=cT$ for
a suitable constant $c\geq 1$ (and we assume $T > c$ in order to have $e^{L/2}<T$). With this choice for
$L$ the quantity $2e^{-L/2}\log \Delta$ becomes $2\log \Delta/\sqrt{cT}$ which is $\geq 1/\sqrt{c}$ under
the assumption that $T \leq 4\log^2\Delta$. Thus, under this further hypothesis on $T$ and with this
choice for $L$, test~\eqref{eq:13A} becomes:
\begin{align}
cT
>&  \sqrt{cT}(2\log \Delta + 2  - 1/\sqrt{c} - \alpha(cT)\,r_1 - \beta(cT)\,n + \log(cT))       \notag \\
&  + 4\sum_{\Norm\a\leq \sqrt{cT}} \Lambda_{\K}(\a)\Big(1 - \frac{\sqrt{cT}}{\Norm\a}\Big)
   + 2\sum_{\substack{T< \Norm\p \leq cT}} \Lambda_{\K}(\p)\log(cT/\Norm\p)
\quad
\implies T_\K \leq T.                                                                           \label{eq:14A}
\end{align}
To appreciate the introduction in this method of the new parameter $c$ we notice that if in~\eqref{eq:14A}
we could retain only the first term appearing on the right hand side then we would conclude that $T_\K
\leq \frac{4}{c} \log^2\Delta$, with a saving with respect to $4$ for every $c > 1$. The two sums on
prime ideals obviously disturb this simple picture, and the main challenge is to get a bound for them
which does not prevent this approach. This is what we do in next section.

\section{Bound for the sums}\label{sec:3A}
The structure of the inequality shows that we should try to estimate the sum of prime ideals via a term
of order $(c-1)^2T$ plus a term depending on the discriminant by at most a quantity of order $(c-1)^2
\sqrt{T} \log\Delta$ or $a(c-1)\sqrt{T}\log\Delta$ for some very small absolute constant $a$: in fact,
if this is the case, these dependencies on $c$ allows to apply the previous approach, at least in the
limit of $c\to 1^+$.\\
The identity
\[
\sum_{\substack{T< \Norm\p \leq cT}} \Lambda_{\K}(\p)\log(cT/\Norm\p)
= \int_T^{cT} \big(\vartheta_\K(u) - \vartheta_\K(T)\big)\frac{\dd u}{u}
\]
where $\vartheta_\K(x):=\sum_{\Norm\p\leq x}\log(\Norm\p)$ and the prime ideal theorem show that the main
term of the sum is $(c-1-\log c)T$, and so it agrees with our assumption. However, we have not been able
to bound the error term better than $a(c-1)\sqrt{T}\log\Delta$ with $a\approx 4/\pi$. The fact that the
bound decays with $c-1$ only to the first order and the value of $a > 1/4$ make this bound completely
useless for our purpose. %
Thus, we change our strategy, and we do not look for the best asymptotic formula, but only for a
convenient upper bound. For this purpose a Brun--Titchmarsh inequality for number fields would be
sufficient, and would produce a result of type $T_\K\leq (4-\alpha)\log^2\Delta$ for some positive
$\alpha$ independent of $n$, but the known results are too coarse for this application (see for
example~\cite{ThornerZaman2}). Thus, we shift to an elementary strategy estimating the sum with the
analogous sum for the rational field: this produces the upper bound we are looking for, but at the cost
of introducing the degree as a factor of the bound.
\begin{proposition}\label{prop:1A}
For every $c\in[1,9/8]$ and $T\geq 73.2$ one has
\[
\frac{1}{n}
\sum_{\substack{T< \Norm\p \leq cT}} \Lambda_{\K}(\p)\log(cT/\Norm\p)
\leq (c-1-\log c)T + \frac{\log c}{8\pi}\sqrt{T}\log^2T + d(c)\,\sqrt{T}.
\]
with $d(c)$ equal to $0.07$, $0.07$, $0.07$ and $0.06$ for $c=9/8$, $13/12$, $17/16$ and $21/20$ (i.e.
for $c=1+1/(4n)$ and $n=2$, $3$, $4$ and $5$) respectively, and $d(c) \leq 0.1$ in any other case.
\end{proposition}
\begin{proof}
For every prime power $p^f$ there are at most $n/f$ prime ideals $\p$ with $\Norm\p = p^f$
%
and for each of them $\Lambda_{\K}(\p)\log(cT/\Norm\p) = f\Lambda(p^f)\log(cT/p^f)$, thus
\[
\sum_{\substack{T< \Norm\p \leq cT}} \Lambda_{\K}(\p)\log(cT/\Norm\p)
\leq n\sum_{T<a\leq cT} \Lambda(a)\log(cT/a).
\]
Integrating by parts, we get
\[
\sum_{\substack{T< a \leq cT}} \Lambda(a)\log(cT/a)
= \int_T^{cT}(\psi(x)-\psi(T))\,\frac{\!\dd x}{x}
= \int_T^{cT}\psi(x)\,\frac{\!\dd x}{x} - \psi(T)\log c,
\]
where $\psi(x) := \sum_{a\leq x}\Lambda(a)$.
To compute the first integral we apply the explicit formula
\[
\psi(x)
= x
  - \lim_{M\to\infty}\sum_{\substack{\rho:\\ |\Imm \rho|\leq M}} \frac{x^\rho}{\rho}
  - \log(2\pi)
  - \frac{1}{2}\log(1-x^{-2})
\]
which is valid whenever $x$ is not a prime power (see~\cite[Ch.~17]{Davenport}). The sum of the last two
terms is negative and the limit converges uniformly on compact sets, therefore an integration term by
term gives
\[
\int_T^{cT}\psi(x)\,\frac{\!\dd x}{x}
\leq
(c-1)T - \sum_{\rho} \frac{(cT)^{\rho}-T^{\rho}}{\rho^2}.
\]
Under RH we can estimate $\sum_{\rho} \frac{(cT)^{\rho}-T^{\rho}}{\rho^2}$ with $\sqrt{T}\sum_{\rho}
\frac{|c^\rho-1|}{|\rho|^2}$. Splitting the sum into two terms according to a parameter $P$, and using
the bound $|c^\rho-1|\leq \sqrt{c}+1$ only for the ones with $|\gamma|>P$ we get
\begin{align*}
T^{-1/2}\Big|\sum_{\rho} \frac{(cT)^\rho-T^\rho}{\rho^2}\Big|
&\leq \sum_{|\gamma|<P} \frac{|c^\rho-1|}{|\rho|^2} + (\sqrt{c}+1)\sum_{|\gamma|>P} \frac{1}{|\rho|^2}                                          \\
&=   -\sum_{|\gamma|<P} \frac{\sqrt{c}+1-|c^\rho-1|}{|\rho|^2} + (\sqrt{c}+1)\sum_{|\gamma|} \frac{1}{|\rho|^2}.
\end{align*}
Under RH the series is equal to $2 - \log(4\pi) + \gamma = 0.04619\ldots$ (compare Equations~(10) and
(11) in~\cite[Ch.~12]{Davenport}),
%
%
%
so that
\[
\int_T^{cT}\psi(x)\,\frac{\!\dd x}{x}
\leq
(c-1)T - \sqrt{T}\sum_{|\gamma|<P} \frac{\sqrt{c}+1-|c^\rho-1|}{|\rho|^2} + 0.0462(\sqrt{c}+1)\sqrt{T}.
\]
Under the assumption that $T\geq 73.2$ (and RH) we also know that
\[
|\psi(T)-T|\leq \frac{\sqrt{T}\log^2 T}{8\pi}
\]
by~\cite[Th.~10]{Schoenfeld1}, so that
\begin{align*}
\int_T^{cT}(\psi(x)-\psi(T))\,\frac{\!\dd x}{x}
\leq& (c-1 - \log c)T
      + \frac{\log c}{8\pi}\sqrt{T}\log^2T                               \\
&     - \sqrt{T}\sum_{|\gamma|<P} \frac{\sqrt{c}+1-|c^\rho-1|}{|\rho|^2}
      + 0.0462(\sqrt{c}+1)\sqrt{T}.
\end{align*}
For a generic $c \leq 9/8$ we pick $P=0$ so that the third sum is zero and $0.0462(\sqrt{c}+1) \leq 0.1$.
%
However, for $c=9/8$, $13/12$, $17/16$ and $21/20$ we pick $P=1500$ and using the known list of zeros for
the Riemann zeta function (for example from~\cite{OdlyzkoTablesZeros}) we get that the $0.1$ constant
improves to $0.07, 0.07, 0.07$ and $0.06$, respectively.
\end{proof}
We also introduce a simple bound for the other sum.
\begin{proposition}\label{prop:2A}
We have
\[
 4\sum_{\Norm\a\leq \sqrt{cT}} \Lambda_{\K}(\a)\Big(1 - \frac{\sqrt{cT}}{\Norm\a}\Big)
 \leq -a_n \sqrt{cT} + b_n
\]
where $(a_2,b_2) = (3.89,82)$, $(a_3,b_3) = (1.93,123)$, $(a_4,b_4) = (1.01,108)$, $(a_5,b_5) = (0.55,135)$
and $(a_n,b_n) = (0,0)$ for every $n\geq 6$.
\end{proposition}
\begin{proof}
For any real $x$, let $x_-:=\min(0,x)$. %
The function $\Lambda_{\K}$ selects ideals $\a$ which are powers of prime ideals $\p^m$, therefore
\begin{align*}
4\sum_{\Norm\a\leq \sqrt{cT}} \Lambda_{\K}(\a)\Big(1 - \frac{\sqrt{cT}}{\Norm\a}\Big)
&  =  4\sum_{\p} \log(\Norm\p)\sum_{m\geq 1}\Big(1 - \frac{\sqrt{cT}}{\Norm\p^m}\Big)_-.
\end{align*}
For every integer prime $p$ there is at least one prime ideal $\p$ sitting above $p$, and its norm is
$p^f$ for some $f\leq n$. Moreover, in all cases $f\big(1 - \frac{\sqrt{cT}}{p^{m f}}\big)_- \leq n\big(1
- \frac{\sqrt{cT}}{p^{m n}}\big)_-$:
to prove it, it is sufficient to prove that $\ell(1 - A B^{-\ell})$ (without the minus) increases along
the positive integers $\ell$ when $A\geq 0$ and $B\geq 2$, i.e. that $\ell(1-A B^{-\ell}) \leq
(\ell+1)(1-A B^{-(\ell+1)})$. This is equivalent to $(1 + \ell - \ell B) A \leq B^{\ell+1}$, which is
true because $1 + \ell \leq \ell B$ in the given ranges for $B$ and $\ell$. Hence
\begin{align*}
4\sum_{\Norm\a\leq \sqrt{cT}} \Lambda_{\K}(\a)\Big(1 - \frac{\sqrt{cT}}{\Norm\a}\Big)
&\leq 4 n\sum_p \log p\sum_{m\geq 1}\Big(1 - \frac{\sqrt{cT}}{p^{n m}}\Big)_-,
\end{align*}
so that
\begin{align*}
4\sum_{\Norm\a\leq \sqrt{cT}} \Lambda_{\K}(\a)\Big(1 - \frac{\sqrt{cT}}{\Norm\a}\Big)
&\leq 4 n\sum_{\substack{p,m\\p^m\in R}}\log p \Big(1 - \frac{\sqrt{cT}}{p^{n m}}\Big)
\end{align*}
for any set of integers $R$. The result follows selecting $R=[1,11]$ for $n=2$ and $3$, $R=[1,8]$ for
$n=4$ and $5$ and $R=\emptyset$ for $n\geq 6$.
%
\end{proof}

\section{Proof of Theorems~\ref{th:1A} and~\ref{th:2A}}\label{sec:4A}
Let $c\in[1,9/8]$ and $T\geq 73.2$. With the bounds in Propositions~\ref{prop:1A} and~\ref{prop:2A},
Test~\eqref{eq:14A} simplifies, after a division by $\sqrt{T}$, to
\begin{align}
c \sqrt{T}
>&  \sqrt{c}(2\log \Delta + 2 - 1/\sqrt{c} - a_n - \alpha(cT)\,r_1 - \beta(cT)\,n + \log(cT))  \notag     \\
&  {+} 2n (c-1-\log c)\sqrt{T} + n\log c\log^2(cT)/(4\pi) + 2d(c)\,n + b_n/\sqrt{T}
\quad
\implies T_\K \leq T.                                                                          \label{eq:15A}
\end{align}
When $n$ and $c$ are fixed and $T$ diverges this implies that
\[
\sqrt{T}\big[c - 2n (c-1-\log c)\big]
>  2\sqrt{c} \log \Delta
   + O_{c,n}(\log^2 T)
\quad
\implies T_\K \leq T,
\]
because $\alpha$ and $\beta$ are bounded. In order to produce a small coefficient for $\log\Delta$ we
need to find $c$ giving a minimum for
\[
f(c,n):=\frac{2\sqrt{c}}{c - 2n (c-1-\log c)}.
\]
When $c$ is close enough to $1$ this quantity is strictly smaller than $2$, for every $n$.
The minimum is attained at a point which is very close to $1+1/(4n)$. Thus we make this choice: $c =
1+1/(4n)$, producing the bound
\[
\sqrt{T}
>  f(1+1/(4n),n) \log \Delta
   + O_{n}(\log^2 T)
\quad
\implies T_\K \leq T.
\]
Elementary tools prove that
\[
f(1+1/(4n),n)^2 < 4-\frac{1}{2n}
\quad\forall n,
\]
%
%
and this proves the existence of a $\Delta_+$ such that~\eqref{eq:1A} holds true for $\Delta \geq
\Delta_+$.
\smallskip\\
For the proof of the range of validity of~\eqref{eq:1A}, we firstly notice that~\cite[Theorem~3.5, (3.10)]
{GrenieMolteni7} with $T_0 = \log\Delta + 2$ proves~\eqref{eq:1A} for any field with $n\leq 8$ and
$\log\Delta\leq 17$.\\
%
We therefore assume henceforth that $\log\Delta\geq 17$, $c=1+1/(4n)$ and $T=(4-1/(2n))\log^2\Delta$.
For degrees $n=2,3,4$ we see that the inequality in~\eqref{eq:15A} holds for $e^{17}\leq \Delta \leq
\Delta_-$, but also for $\Delta \geq \Delta_+$ with $\Delta_\pm$ as given in Table~\ref{tab:1A}.
With the previous computation, this completes the proof of the claim in Theorem~\ref{th:1A} for these
degrees. %
\smallskip\\
For $n= 5,6,7$ and $8$, \eqref{eq:15A} holds for all $\log \Delta \geq 17$ without exceptions,
%
%
%
%
%
hence the claim in Theorem~\ref{th:1A} is proved also for these degrees. %
\smallskip\\
We now prove that~\eqref{eq:1A} holds when $n\geq 9$: this will complete the proof of Theorem~\ref{th:1A}.\\
From Odlyzko computations~\cite{Odlyzko:tables,Odlyzko2,Odlyzko3} any field of degree $n\geq 9$ satisfies
$\log\Delta\geq 17$
%
%
so that for $c = 1+1/(4n)$ and $T=(4-1/(2n))\log^2\Delta$ one has that $cT \geq 4\log^2\Delta \geq 4\cdot
17^2 = 34^2$. Hence $\alpha(cT)\geq \alpha(34^2) \geq 1$ and $\beta(cT)\geq \beta(34^2) \geq 4.42$,
so that~\eqref{eq:15A} simplifies to
\begin{align}
c \sqrt{T}
>&  \sqrt{c}(2\log \Delta + 2 - 1/\sqrt{c} - a_n - r_1 - 4.42\,n + \log(cT))              \notag     \\
&  {+} 2n (c-1-\log c)\sqrt{T} + n\log c\log^2(cT)/(4\pi) + 2d(c)\,n + b_n/\sqrt{T}
\quad
\implies T_\K \leq T.                                                                     \label{eq:16A}
\end{align}
Set $S := \sqrt{cT} = [(1+1/(4n))(4-1/(2n))]^{1/2}\log\Delta$.
We rewrite~\eqref{eq:16A} in terms of $S$; then we simplify a bit the resulting inequality dividing by
$\sqrt{c}$, noticing that $n\log c\leq n(c-1)\leq 1/4$, using the bound $1/\sqrt{c}\geq
1/\sqrt{1+1/(4\cdot 9)}\geq 0.98$ and removing the $r_1$ term. This produces the test:
\begin{multline}\label{eq:17A}
\Big[1 - 2n\frac{c-1-\log c}{c} - \frac{2/\sqrt{c}}{\sqrt{4-1/(2n)}}\Big] S\\
> 1.02 - 4.42\,n + \frac{2d(c)}{\sqrt{c}}\,n + 2\log S + \frac{1/\sqrt{c}}{4\pi}\log^2 S
\quad
\implies T_\K \leq T.
\end{multline}
Finally, we further simplify the inequality noticing that the coefficient of $S$ on the left hand side is
larger than $4/(10^3 n^2)$ when $n\geq 9$ and removing the $1/\sqrt{c}\leq 1$ factor from the right.
%
%
This yields:
\begin{equation}\label{eq:18A}
\frac{4/10^3}{n^2} S + 4.22\,n
> 1.02 + 2\log S + \frac{1}{4\pi}\log^2S
\quad
\implies T_\K \leq T.
\end{equation}
%
%
In terms of $n$, the function $An^{-2}+Bn$ has a minimum at $n = n_0 := (2A/B)^{1/3}$, with a value
which is $3(AB^2/4)^{1/3}$. Thus, \eqref{eq:18A} is true as soon as
\[
0.48 S > (1.02 + 2\log S + (\log^2S)/(4\pi))^3,
\]
%
and for this it is sufficient to assume that $S \geq 80100$.
%
Suppose $S < 80100$. Then $n_0 = (2A/B)^{1/3}$ is $\leq 6$.
However we are assuming that $n\geq 9$, hence in our setting ($S < 80100$ and $n\geq 9$) the minimum of
the function appearing on the left hand side of~\eqref{eq:18A} is not attained at $n=n_0$ but at $n=9$.
The inequality becomes
\[
\frac{4/10^3}{9^2} S + 4.22\cdot 9
> 1.02 + 2\log S + \frac{1}{4\pi}\log^2S,
\]
which is true.
\medskip\\
The proof of Theorem~\ref{th:2A} is similar and even simpler. By Theorem~\ref{th:1A} the claim in
Theorem~\ref{th:2A} is true for if $n\geq 5$ or, for $2\leq n\leq 4$, if $\Delta \leq \Delta_-$, hence we
can assume that $\log\Delta \geq 63$. Let $c = 1+1/(4n)$ and $T=(4-1/(2n)+1/(2n^2))\log^2\Delta$. The
assumption $cT \geq 34^2$ is satisfied, hence we have~\eqref{eq:16A} at our disposal. As we have done
for Theorem~\ref{th:1A}, we rewrite~\eqref{eq:16A} in terms of $S := \sqrt{cT} =
[(1+1/(4n))(4-1/(2n)+1/(2n^2))]^{1/2} \log\Delta$, which is $\geq 126$, because $S\geq 2\log\Delta$. With
the same steps as for~\eqref{eq:17A} this yields
%
%
\begin{multline*}
\Big[1 - 2n\frac{c-1-\log c}{c} - \frac{2/\sqrt{c}}{(4{-}\frac{1}{2n}{+}\frac{1}{2n^2})^{1/2}}\Big] S\\
>  \frac{b_n}{S} + 2 - \frac{1}{\sqrt{c}} - a_n - 4.22\,n + 2\log S + \frac{1}{4\pi}\log^2S.
\end{multline*}
The coefficient of $S$ on the left hand side is $\geq 0.06/n^2$ for $n=2,3,4,5$
%
%
and $b_n/S + 2 - 1/\sqrt{c} \leq 2.1$ for the same degrees, hence in this case the test simplifies to
%
\[
\frac{0.06}{n^2} S
> 2.1 - a_n - 4.22\,n + 2\log S + \frac{1}{4\pi}\log^2 S
\quad
\implies T_\K \leq T.
\]
This inequality is indeed true for $S\geq 126$ and $n=2,3,4$.
%

\section{Proof of Theorems~\ref{th:3A} and~\ref{th:4A}}\label{sec:5A}
We recall that in any finite abelian group $G$ every subgroup $H$ can be realized as intersection of
kernels of characters in a suitable set $S=S(H)$, and that the mapping $H \mapsto S(H)$ is
inclusion-reversing.
In this correspondence $\Cl^2$ is mapped to $\{\chi\colon \chi^2=\chi_0\}$.
Thus, suppose that $T < T'_\K$, so that the prime ideals with norm $\leq T$ generate a subgroup $H$ which
does not contain $\Cl^2$. Then the set of characters which are trivial on $H$ is not contained into the
set of characters which are trivial on $\Cl^2$, i.e. in the set of characters whose square is trivial.
Hence there exists a character $\chi$ with $\chi^2 \neq \chi_0$ and which is trivial on $H$. Thus, both
$\chi$ and $\chi^2$ are not trivial but are trivial on $H$ so that we have at our disposal the
relation~\eqref{eq:4A} for both of them, under the same assumption for $L$ and $T$. In particular, we
combine linearly the two relations with coefficient $2/3$ for the one for $\chi$ and $1/3$ for the one
for $\chi^2$. This combination produces the equality
\begin{align*}
\phi(1) + \phi(0) - \sum_{\rho_\K} \phi(\rho_\K)
&+ \frac{2}{3}\sum_{\rho_\chi} \phi(\rho_\chi)
 + \frac{1}{3}\sum_{\rho_{\chi^2}} \phi(\rho_{\chi^2}) \\
&= \sum_{\substack{T< \Norm\p \leq e^L}} \frac{\Lambda_{\K}(\p)}{\sqrt{\Norm\p}}\frac{2|1-\chi(\p)|^2 + |1-\chi^2(\p)|^2 }{3} F(\log \Norm\p).
\end{align*}
The maximum of $2|1-z|^2 + |1-z^2|^2$ for $|z|=1$ is $9$.
%
Thus, proceeding as we have done for the proof of Theorem~\ref{th:1A}, we conclude that in case
\[
4\int_0^{+\infty}F(x)\Ch(x/2)\dd x
> \sum_{\rho_\K} \phi(\rho_\K)
    + 3\sum_{\substack{T< \Norm\p \leq e^L}} \frac{\Lambda_{\K}(\p)}{\sqrt{\Norm\p}}F(\log \Norm\p)
\quad \text{then} \quad T'_\K \leq T.
\]
The difference with respect to~\eqref{eq:6A} is due to the fact that now in front of the sum on prime
ideals we have $3$ in place of $4$. Keeping the same steps and with the same choices for $\psi^+$, $L$,
$c$ and $T$ we get that~\eqref{eq:15A} changes into
\begin{align}
c \sqrt{T}
>&  \sqrt{c}(2\log \Delta + 2 - 1/\sqrt{c} - a_n - \alpha(cT)\,r_1 - \beta(cT)\,n + \log(cT))  \notag     \\
&  {+} \frac{3}{2}\,n (c-1-\log c)\sqrt{T} + \frac{3}{16\pi}\,n\log c\log^2(cT) + \frac{3d(c)}{2}\,n + b_n/\sqrt{T}
\quad
\implies T'_\K \leq T.                                                                         \label{eq:19A}
\end{align}
Now we proceed as for Theorems~\ref{th:1A} and~\ref{th:2A}, but starting with~\eqref{eq:19A}. For
example, to produce a small coefficient for $\log\Delta$ we need to find $c$ giving a minimum for
\[
f(c,n):=\frac{2\sqrt{c}}{c - \frac{3}{2}\,n (c-1-\log c)}.
\]
For $c = 1+1/(3n)$ this gives the bound
\[
\sqrt{T}
>  f(1+1/(3n),n) \log \Delta
   + O_{n}(\log^2 T)
\quad
\implies T'_\K \leq T.
\]
Since
\[
f(1+1/(3n),n)^2 < 4-\frac{2}{3n}
\quad\forall n,
\]
this proves the existence of a $\Delta'_{+}$ such that~\eqref{eq:2A} holds true for $\Delta \geq
\Delta'_{+}$. All other claims in Theorem~\ref{th:3A} are proved imitating what we have done for
Theorem~\ref{th:1A}, and Theorem~\ref{th:4A} is proved imitating what we have done for
Theorem~\ref{th:2A}.

\section{Algorithmic version}\label{sec:6A}
The algorithm is based on the following result. It is a version of the main result, Theorem~2.1,
of~\cite{BelabasDiazyDiazFriedman}, once one decides to allow the support of the test function to be
larger than $\log T$.
\begin{theorem}\label{th:5A}
Let $\K$ be a number field satisfying the Riemann Hypothesis for all $L$-functions attached to
non-trivial characters of its ideal class group $\Cl$. Let $E$ be a set of prime ideals of $\K$ and let
$X_E$ be its characteristic function. Suppose that there exists a non-negative function $F\in\W$ of
compact support with non-negative Fourier transform such that
\begin{equation}
2\sum_{\p}(-1)^{X_E(\p)}\log\Norm\p\sum_{m\geq 1}\frac{F(m\log\Norm\p)}{\sqrt{\Norm\p^m}}
> F(0)(\log\Delta - (\gamma + \log 8\pi)n)
+ I(F) n
- J(F) r_1.                                                                                  \label{eq:20A}
\end{equation}
Then $\{\p\colon \Norm\p\in\supp F\}\backslash E$ is a generating set for $\Cl$.
\end{theorem}
Above and below, all sums are intended limited to the primes in the support of $F$ (otherwise the summand
is $0$).
\begin{proof}
Let $E$ be a set of primes of $\K$, $X_E$ be its characteristic function and $F\in\W$ be a non-negative
function with compact support and non-negative Fourier transform. Suppose that $G:=\{\p\colon
\Norm\p\in\supp F\}\backslash E$ does not generate $\Cl$. Then it generates a proper subgroup of $\Cl$.
Hence there exists a character $\chi\neq \chi_0$ of $\Cl$ which is trivial on said subgroup. Then
from~\eqref{eq:3A} we have
\begin{multline*}
\sum_{\a} \frac{\Lambda_{\K}(\a)}{\sqrt{\Norm\a}}(\chi(\a)+\overline{\chi(\a)}) F(\log \Norm\a) \\
= F(0)(\log \Delta - (\gamma+\log 8\pi)n)
  + I(F) n
  - J(F) r_1
  - \sum_{\rho_\chi} \phi(\rho_\chi).
\end{multline*}
We have $\Lambda_\K(\a)=0$ unless $\a$ is the power of a prime ideal $\p$, and $\chi(\p)=1$ if
$\p\in G$. Therefore
\begin{multline*}
2\sum_{\p\in G} \log\Norm\p\sum_{m\geq 1}\frac{F(m\log \Norm\p)}{\sqrt{\Norm\p^m}}
+\sum_{\p\in E} \log\Norm\p\sum_{m\geq 1}\frac{(\chi(\p^m)+\overline{\chi(\p^m)})F(\log m\Norm\p)}{\sqrt{\Norm\p^m}} \\
= F(0)(\log \Delta - (\gamma+\log 8\pi)n)
  + I(F) n
  - J(F) r_1
  - \sum_{\rho_\chi} \phi(\rho_\chi).
\end{multline*}
Now, since we assumed GRH and that the Fourier transform of $F$ is non negative, we have
$\phi(\rho_\chi)\geq 0$ and, moreover, $\chi(\p^m)+\overline{\chi(\p^m)}\geq -2$, therefore
\begin{align*}
&  2\sum_{\p} (-1)^{X_E(\p)}\log\Norm\p\sum_{m\geq 1}\frac{F(m\log \Norm\p)}{\sqrt{\Norm\p^m}} \\
&= 2\sum_{\p\in G} \log\Norm\p\sum_{m\geq 1}\frac{F(m\log \Norm\p)}{\sqrt{\Norm\p^m}}
  -2\sum_{\p\in E} \log\Norm\p\sum_{m\geq 1}\frac{F(m\log \Norm\p)}{\sqrt{\Norm\p^m}}          \\
&\leq F(0)(\log \Delta - (\gamma+\log 8\pi)n)
  + I(F) n
  - J(F) r_1,
\end{align*}
where we have used the fact that $\p\in G$ if and only if $\p\in\supp F$ and $X_E(\p)=0$. Therefore, if
the opposite inequality is satisfied, there does not exist any such $\chi$, hence $G$ generates the full
class group.
\end{proof}
To use this theorem, the idea is to take a function with large enough support so that the sum on prime
ideals out of $E$ minus the right hand side of~\eqref{eq:20A} is positive. One can then ``use'' that
positive number to ``remove'' the ideals in $E$ from the left hand side. Since we are interested in
reducing the bound for the norm of the ideals in the generating set, we will take $E$ of the form
\[
E:=\{\p\colon\Norm\p> T\}.
\]
To simplify exposition (and since we do not really need the other cases anyway!), we will suppose that
$\supp F=[-\log cT,\log cT]$ with $1\leq c<\sqrt{T}$, hence when $X_E(\p)=1$, the only value of $m$ that
matters is $m=1$.

%
Suppose that $F$ is as in the hypothesis of the theorem and that~\eqref{eq:20A} is satisfied, and $E$ and
$cT$ are as explained. Subtracting~\eqref{eq:3A} for $\chi=\chi_0$ from~\eqref{eq:20A}, one finds that
\[
4\int_0^{+\infty}F(x)\Ch(x/2)\dd x=\phi(0)+\phi(1)>\sum_{\rho_\K}\phi(\rho_\K)
+ 4\sum_{T<\Norm\p\leq cT}\log\Norm\p\frac{F(\log\Norm\p)}{\sqrt{\Norm\p}},
\]
which is the opposite version of~\eqref{eq:5A}. To find a good test function $F$, one can then follow the
same path as the one of the proof of Theorem~\ref{th:1A}. The conclusion is then that such a good
function is exactly the one that has been chosen in~\eqref{eq:10A}. The corresponding result is as
follows.
\begin{corollary}\label{cor:1A}
Let $T\geq 1$, $c\in[1,T)$ and $L=\log(cT)$. Suppose
\begin{align*}
\sum_{\Norm\p^m\leq \sqrt{cT}}\Bigl(m\log\Norm\p-2+2\frac{\sqrt{cT}}{\Norm\p^m}\Bigr)\log\Norm\p
&+  \sum_{\sqrt{cT}< \Norm\p^m\leq cT}\log(cT/\Norm\p^m)\log\Norm\p                               \\
&- 2\sum_{T< \Norm\p\leq cT}\log(cT/\Norm\p)\log\Norm\p                                           \\
>
(\sqrt{cT}-1)
\bigl(\log\Delta-\bigl(\gamma+\log\bigl(&2\pi\bigl(1-(cT)^{-1/2}\bigr)\bigr)\bigr)n - r_1\log 2\bigr)
- \frac{L(L-4)}{8}(n+r_1).
\end{align*}
Then $T_\K\leq T$.
\end{corollary}
\begin{proof}
Let $T\geq 1$, $c\in[1,T)$, $L=\log(cT)$ and let $F\in \W$, supported in $[-\log cT, \log cT]$, be as
in~\eqref{eq:10A}. From Theorem~\ref{th:5A}, applied to $F$, we deduce that if
\begin{align*}
2\sum_{\Norm\p^m\leq \sqrt{cT}}\Bigl(m\log\Norm\p - 2 + 2\frac{\sqrt{cT}}{\Norm\p^m}\Bigr)\log\Norm\p
&+ 2\sum_{\sqrt{cT}< \Norm\p^m\leq cT}\log(cT/\Norm\p^m)\log\Norm\p                                \\
&- 4\sum_{T< \Norm\p\leq cT}\log(cT/\Norm\p)\log\Norm\p                                            \\
>
2(\sqrt{cT}-1)(\log\Delta - (\gamma + \log &8\pi)n)
+ I(F) n
- J(F) r_1,
\end{align*}
then $T_\K\leq T$. From~\eqref{eq:11aA}--\eqref{eq:11bA} we have
\[
I(F) n - J(F) r_1
\leq 2(\sqrt{cT}-1)(\log(4/(1-(cT)^{-1/2}))n - r_1\log 2)-L(L-4)(n+r_1)/4.
\]
%
%
%
%
The result follows.
\end{proof}
Thus, we keep $F$ as in~\eqref{eq:10A} and set
\begin{align*}
D_\K(cT)
&:=
\frac{\sqrt{cT}-1}{2}\big(\log\Delta - \bigl(\gamma + \log\bigl(2\pi\bigl(1-(cT)^{-1/2}\bigr)\bigr)\bigr)n - r_1 \log 2\big) \\
&\qquad
- \frac{L(L-4)}{8}(n+r_1)
- \frac{1}{2}\sum_{\p}\log\Norm\p\sum_{m\geq 1}\frac{F(m\log\Norm\p)}{\sqrt{\Norm\p^m}}
\intertext{and}
S_\K(T,cT)
&:=
\sum_{T< \Norm\p\leq cT}\log(cT/\Norm\p)\log\Norm\p.
\end{align*}
Observe that $D_\K$ depends on $cT$ but not on $T$. The inequality in the hypothesis of the corollary is
satisfied if $D_\K(cT)+S_\K(T,cT)<0$. \\
The algorithm runs as follows: firstly it determines an integer $x$ such that $D_\K(x) < 0$. Then it
computes the sum $S_\K(t,x)$ with $t$ starting from $x$ and decreasing until $D_\K(x)+S_\K(t,x)$ turns
positive. Then $T:=t+1$ is the upper bound for $T_\K$ and $c = x/T$.
The resulting $T$ depends on $x$, so that the algorithm chooses an $x$ giving the lowest $T$.
\medskip

To ease the comparison of the new algorithm with the previous ones we introduce some notation. So, let
$T(\K)$ be the result of the original algorithm of~\cite{BelabasDiazyDiazFriedman}, $T_1(\K)$ the result
of the optimal algorithm of~\cite{GrenieMolteni7}, $T_2(\K)$ the result of the simplified algorithm of
the same paper and $T_{2+}(\K)$ the result of the algorithm described above.\\
From the discussion between Theorem~\ref{th:5A} and Corollary~\ref{cor:1A}, we have immediately the
following upper bound for $T_{2+}(\K)$.
\begin{corollary}\label{cor:2A}
The bounds for $T_\K$ given in Theorems~\ref{th:1A} and~\ref{th:2A} are valid also for $T_{2+}(\K)$.
\end{corollary}
We have tested the algorithms for several families of fields, getting in all cases results which are very
similar to the ones we discuss for the three series we reproduce in details in next figures, which are:
the quadratic fields of the form $\Q(x)/(x^2-p_n)$, the sextic fields which are the Galois
closures of $\Q(x)/(x^3-p_n)$ and fields in degree 21 of the form $\Q(x)/(x^{21}-p_n)$; in
all cases $p_n$ is the first strong pseudo prime after $10^n$ and all series start with $n=0$.\\
All computations have been performed on the PlaFRIM cluster~\cite{BORDEAUX-CLUSTER}.

\enlargethispage{2\baselineskip}
\begin{fixedfig}
\begingroup
  \makeatletter
  \providecommand\color[2][]{%
    \GenericError{(gnuplot) \space\space\space\@spaces}{%
      Package color not loaded in conjunction with
      terminal option `colourtext'%
    }{See the gnuplot documentation for explanation.%
    }{Either use 'blacktext' in gnuplot or load the package
      color.sty in LaTeX.}%
    \renewcommand\color[2][]{}%
  }%
  \providecommand\includegraphics[2][]{%
    \GenericError{(gnuplot) \space\space\space\@spaces}{%
      Package graphicx or graphics not loaded%
    }{See the gnuplot documentation for explanation.%
    }{The gnuplot epslatex terminal needs graphicx.sty or graphics.sty.}%
    \renewcommand\includegraphics[2][]{}%
  }%
  \providecommand\rotatebox[2]{#2}%
  \@ifundefined{ifGPcolor}{%
    \newif\ifGPcolor
    \GPcolorfalse
  }{}%
  \@ifundefined{ifGPblacktext}{%
    \newif\ifGPblacktext
    \GPblacktexttrue
  }{}%
  \let\gplgaddtomacro\g@addto@macro
  \gdef\gplbacktext{}%
  \gdef\gplfronttext{}%
  \makeatother
  \ifGPblacktext
    \def\colorrgb#1{}%
    \def\colorgray#1{}%
  \else
    \ifGPcolor
      \def\colorrgb#1{\color[rgb]{#1}}%
      \def\colorgray#1{\color[gray]{#1}}%
      \expandafter\def\csname LTw\endcsname{\color{white}}%
      \expandafter\def\csname LTb\endcsname{\color{black}}%
      \expandafter\def\csname LTa\endcsname{\color{black}}%
      \expandafter\def\csname LT0\endcsname{\color[rgb]{1,0,0}}%
      \expandafter\def\csname LT1\endcsname{\color[rgb]{0,1,0}}%
      \expandafter\def\csname LT2\endcsname{\color[rgb]{0,0,1}}%
      \expandafter\def\csname LT3\endcsname{\color[rgb]{1,0,1}}%
      \expandafter\def\csname LT4\endcsname{\color[rgb]{0,1,1}}%
      \expandafter\def\csname LT5\endcsname{\color[rgb]{1,1,0}}%
      \expandafter\def\csname LT6\endcsname{\color[rgb]{0,0,0}}%
      \expandafter\def\csname LT7\endcsname{\color[rgb]{1,0.3,0}}%
      \expandafter\def\csname LT8\endcsname{\color[rgb]{0.5,0.5,0.5}}%
    \else
      \def\colorrgb#1{\color{black}}%
      \def\colorgray#1{\color[gray]{#1}}%
      \expandafter\def\csname LTw\endcsname{\color{white}}%
      \expandafter\def\csname LTb\endcsname{\color{black}}%
      \expandafter\def\csname LTa\endcsname{\color{black}}%
      \expandafter\def\csname LT0\endcsname{\color{black}}%
      \expandafter\def\csname LT1\endcsname{\color{black}}%
      \expandafter\def\csname LT2\endcsname{\color{black}}%
      \expandafter\def\csname LT3\endcsname{\color{black}}%
      \expandafter\def\csname LT4\endcsname{\color{black}}%
      \expandafter\def\csname LT5\endcsname{\color{black}}%
      \expandafter\def\csname LT6\endcsname{\color{black}}%
      \expandafter\def\csname LT7\endcsname{\color{black}}%
      \expandafter\def\csname LT8\endcsname{\color{black}}%
    \fi
  \fi
    \setlength{\unitlength}{0.0500bp}%
    \ifx\gptboxheight\undefined%
      \newlength{\gptboxheight}%
      \newlength{\gptboxwidth}%
      \newsavebox{\gptboxtext}%
    \fi%
    \setlength{\fboxrule}{0.5pt}%
    \setlength{\fboxsep}{1pt}%
\begin{picture}(7936.00,3400.00)%
    \gplgaddtomacro\gplbacktext{%
      \csname LTb\endcsname%
      \put(462,440){\makebox(0,0)[r]{\strut{}$0$}}%
      \put(462,825){\makebox(0,0)[r]{\strut{}$5$}}%
      \put(462,1210){\makebox(0,0)[r]{\strut{}$10$}}%
      \put(462,1595){\makebox(0,0)[r]{\strut{}$15$}}%
      \put(462,1980){\makebox(0,0)[r]{\strut{}$20$}}%
      \put(462,2365){\makebox(0,0)[r]{\strut{}$25$}}%
      \put(462,2750){\makebox(0,0)[r]{\strut{}$30$}}%
      \put(462,3135){\makebox(0,0)[r]{\strut{}$35$}}%
      \put(657,220){\makebox(0,0){\strut{}$0$}}%
      \put(1640,220){\makebox(0,0){\strut{}$250$}}%
      \put(2623,220){\makebox(0,0){\strut{}$500$}}%
      \put(3606,220){\makebox(0,0){\strut{}$750$}}%
      \put(4590,220){\makebox(0,0){\strut{}$1000$}}%
      \put(5573,220){\makebox(0,0){\strut{}$1250$}}%
      \put(6556,220){\makebox(0,0){\strut{}$1500$}}%
      \put(7539,220){\makebox(0,0){\strut{}$1750$}}%
    }%
    \gplgaddtomacro\gplfronttext{%
      \colorrgb{0.58,0.00,0.83}%
      \put(6219,2962){\makebox(0,0)[l]{\strut{}quadratic}}%
      \colorrgb{0.00,0.62,0.45}%
      \put(6219,2742){\makebox(0,0)[l]{\strut{}degree 6}}%
      \colorrgb{0.34,0.71,0.91}%
      \put(6219,2522){\makebox(0,0)[l]{\strut{}degree 21}}%
    }%
    \gplbacktext
    \put(0,0){\includegraphics{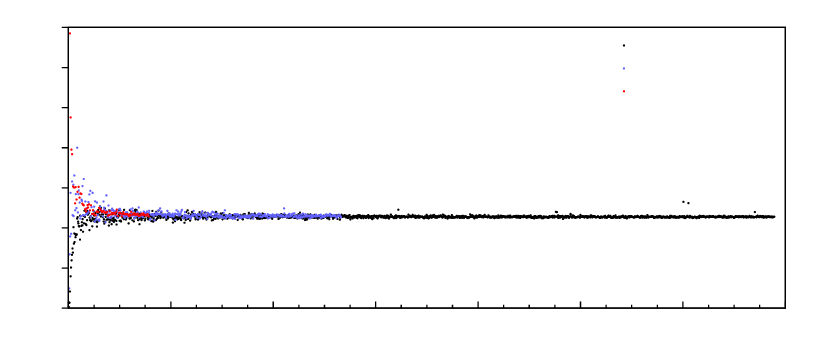}}%
    \gplfronttext
  \end{picture}%
\endgroup
\hspace{1cm}\mbox{.}\\
\caption{$\dfrac{T_{2+}(\K)}{T(\K)}(\log\log\Delta)^2$, in abscissa $\frac{1}{n}\log\Delta$.}%
\label{fig:1A}
\end{fixedfig}
\smallskip
\noindent
Figure~\ref{fig:1A} shows that $T_{2+}(\K)$ improves on $T(\K)$ by a factor which is the square of a
double log of $\Delta$: this fact was already known from~\cite{GrenieMolteni7} for $T_1(\K)$ and
$T_2(\K)$, and $T_{2+}(\K)\leq T_2(\K)$ by design.
\smallskip

\begin{fixedfig}
\begingroup
  \makeatletter
  \providecommand\color[2][]{%
    \GenericError{(gnuplot) \space\space\space\@spaces}{%
      Package color not loaded in conjunction with
      terminal option `colourtext'%
    }{See the gnuplot documentation for explanation.%
    }{Either use 'blacktext' in gnuplot or load the package
      color.sty in LaTeX.}%
    \renewcommand\color[2][]{}%
  }%
  \providecommand\includegraphics[2][]{%
    \GenericError{(gnuplot) \space\space\space\@spaces}{%
      Package graphicx or graphics not loaded%
    }{See the gnuplot documentation for explanation.%
    }{The gnuplot epslatex terminal needs graphicx.sty or graphics.sty.}%
    \renewcommand\includegraphics[2][]{}%
  }%
  \providecommand\rotatebox[2]{#2}%
  \@ifundefined{ifGPcolor}{%
    \newif\ifGPcolor
    \GPcolorfalse
  }{}%
  \@ifundefined{ifGPblacktext}{%
    \newif\ifGPblacktext
    \GPblacktexttrue
  }{}%
  \let\gplgaddtomacro\g@addto@macro
  \gdef\gplbacktext{}%
  \gdef\gplfronttext{}%
  \makeatother
  \ifGPblacktext
    \def\colorrgb#1{}%
    \def\colorgray#1{}%
  \else
    \ifGPcolor
      \def\colorrgb#1{\color[rgb]{#1}}%
      \def\colorgray#1{\color[gray]{#1}}%
      \expandafter\def\csname LTw\endcsname{\color{white}}%
      \expandafter\def\csname LTb\endcsname{\color{black}}%
      \expandafter\def\csname LTa\endcsname{\color{black}}%
      \expandafter\def\csname LT0\endcsname{\color[rgb]{1,0,0}}%
      \expandafter\def\csname LT1\endcsname{\color[rgb]{0,1,0}}%
      \expandafter\def\csname LT2\endcsname{\color[rgb]{0,0,1}}%
      \expandafter\def\csname LT3\endcsname{\color[rgb]{1,0,1}}%
      \expandafter\def\csname LT4\endcsname{\color[rgb]{0,1,1}}%
      \expandafter\def\csname LT5\endcsname{\color[rgb]{1,1,0}}%
      \expandafter\def\csname LT6\endcsname{\color[rgb]{0,0,0}}%
      \expandafter\def\csname LT7\endcsname{\color[rgb]{1,0.3,0}}%
      \expandafter\def\csname LT8\endcsname{\color[rgb]{0.5,0.5,0.5}}%
    \else
      \def\colorrgb#1{\color{black}}%
      \def\colorgray#1{\color[gray]{#1}}%
      \expandafter\def\csname LTw\endcsname{\color{white}}%
      \expandafter\def\csname LTb\endcsname{\color{black}}%
      \expandafter\def\csname LTa\endcsname{\color{black}}%
      \expandafter\def\csname LT0\endcsname{\color{black}}%
      \expandafter\def\csname LT1\endcsname{\color{black}}%
      \expandafter\def\csname LT2\endcsname{\color{black}}%
      \expandafter\def\csname LT3\endcsname{\color{black}}%
      \expandafter\def\csname LT4\endcsname{\color{black}}%
      \expandafter\def\csname LT5\endcsname{\color{black}}%
      \expandafter\def\csname LT6\endcsname{\color{black}}%
      \expandafter\def\csname LT7\endcsname{\color{black}}%
      \expandafter\def\csname LT8\endcsname{\color{black}}%
    \fi
  \fi
    \setlength{\unitlength}{0.0500bp}%
    \ifx\gptboxheight\undefined%
      \newlength{\gptboxheight}%
      \newlength{\gptboxwidth}%
      \newsavebox{\gptboxtext}%
    \fi%
    \setlength{\fboxrule}{0.5pt}%
    \setlength{\fboxsep}{1pt}%
\begin{picture}(7936.00,3400.00)%
    \gplgaddtomacro\gplbacktext{%
      \csname LTb\endcsname%
      \put(594,440){\makebox(0,0)[r]{\strut{}$0$}}%
      \put(594,710){\makebox(0,0)[r]{\strut{}$0.1$}}%
      \put(594,979){\makebox(0,0)[r]{\strut{}$0.2$}}%
      \put(594,1249){\makebox(0,0)[r]{\strut{}$0.3$}}%
      \put(594,1518){\makebox(0,0)[r]{\strut{}$0.4$}}%
      \put(594,1788){\makebox(0,0)[r]{\strut{}$0.5$}}%
      \put(594,2057){\makebox(0,0)[r]{\strut{}$0.6$}}%
      \put(594,2327){\makebox(0,0)[r]{\strut{}$0.7$}}%
      \put(594,2596){\makebox(0,0)[r]{\strut{}$0.8$}}%
      \put(594,2865){\makebox(0,0)[r]{\strut{}$0.9$}}%
      \put(594,3135){\makebox(0,0)[r]{\strut{}$1$}}%
      \put(789,220){\makebox(0,0){\strut{}$0$}}%
      \put(1753,220){\makebox(0,0){\strut{}$250$}}%
      \put(2718,220){\makebox(0,0){\strut{}$500$}}%
      \put(3682,220){\makebox(0,0){\strut{}$750$}}%
      \put(4646,220){\makebox(0,0){\strut{}$1000$}}%
      \put(5610,220){\makebox(0,0){\strut{}$1250$}}%
      \put(6575,220){\makebox(0,0){\strut{}$1500$}}%
      \put(7539,220){\makebox(0,0){\strut{}$1750$}}%
    }%
    \gplgaddtomacro\gplfronttext{%
      \colorrgb{0.58,0.00,0.83}%
      \put(6219,1053){\makebox(0,0)[l]{\strut{}quadratic}}%
      \colorrgb{0.00,0.62,0.45}%
      \put(6219,833){\makebox(0,0)[l]{\strut{}degree 6}}%
      \colorrgb{0.34,0.71,0.91}%
      \put(6219,613){\makebox(0,0)[l]{\strut{}degree 21}}%
    }%
    \gplbacktext
    \put(0,0){\includegraphics{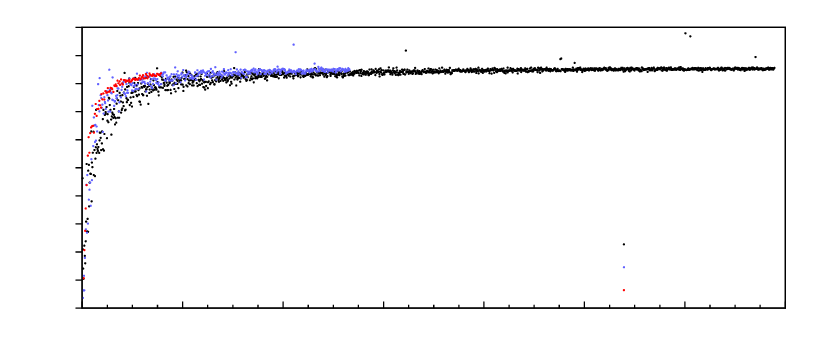}}%
    \gplfronttext
  \end{picture}%
\endgroup
\hspace{1cm}\mbox{.}\\
\caption{$\dfrac{T_{2+}(\K)}{\log^2\Delta}$, in abscissa $\frac{1}{n}\log\Delta$.}
\label{fig:2A}
\end{fixedfig}
\smallskip
\noindent
Figure~\ref{fig:2A} shows that $T_{2+}(\K)$ is significantly smaller than $\log^2\Delta$; this fact
suggests that Theorems~\ref{th:1A}--\ref{th:4A} still do not fully describe the real behavior of $T_\K$,
which is not surprising because in going from~\eqref{eq:7A} to~\eqref{eq:8A}, we major
$|\Ree\widehat{\psi^+}|$ by $|\widehat{\psi^+}|$, which means that $\sqrt{cT}$ is multiplied by $2$ when
$|\Imm\widehat{\psi^+}|=|\Ree\widehat{\psi^+}|$.
\smallskip

\enlargethispage{8\baselineskip}
\begin{fixedfig}
\begingroup
  \inputencoding{cp1252}%
  \makeatletter
  \providecommand\color[2][]{%
    \GenericError{(gnuplot) \space\space\space\@spaces}{%
      Package color not loaded in conjunction with
      terminal option `colourtext'%
    }{See the gnuplot documentation for explanation.%
    }{Either use 'blacktext' in gnuplot or load the package
      color.sty in LaTeX.}%
    \renewcommand\color[2][]{}%
  }%
  \providecommand\includegraphics[2][]{%
    \GenericError{(gnuplot) \space\space\space\@spaces}{%
      Package graphicx or graphics not loaded%
    }{See the gnuplot documentation for explanation.%
    }{The gnuplot epslatex terminal needs graphicx.sty or graphics.sty.}%
    \renewcommand\includegraphics[2][]{}%
  }%
  \providecommand\rotatebox[2]{#2}%
  \@ifundefined{ifGPcolor}{%
    \newif\ifGPcolor
    \GPcolorfalse
  }{}%
  \@ifundefined{ifGPblacktext}{%
    \newif\ifGPblacktext
    \GPblacktexttrue
  }{}%
  \let\gplgaddtomacro\g@addto@macro
  \gdef\gplbacktext{}%
  \gdef\gplfronttext{}%
  \makeatother
  \ifGPblacktext
    \def\colorrgb#1{}%
    \def\colorgray#1{}%
  \else
    \ifGPcolor
      \def\colorrgb#1{\color[rgb]{#1}}%
      \def\colorgray#1{\color[gray]{#1}}%
      \expandafter\def\csname LTw\endcsname{\color{white}}%
      \expandafter\def\csname LTb\endcsname{\color{black}}%
      \expandafter\def\csname LTa\endcsname{\color{black}}%
      \expandafter\def\csname LT0\endcsname{\color[rgb]{1,0,0}}%
      \expandafter\def\csname LT1\endcsname{\color[rgb]{0,1,0}}%
      \expandafter\def\csname LT2\endcsname{\color[rgb]{0,0,1}}%
      \expandafter\def\csname LT3\endcsname{\color[rgb]{1,0,1}}%
      \expandafter\def\csname LT4\endcsname{\color[rgb]{0,1,1}}%
      \expandafter\def\csname LT5\endcsname{\color[rgb]{1,1,0}}%
      \expandafter\def\csname LT6\endcsname{\color[rgb]{0,0,0}}%
      \expandafter\def\csname LT7\endcsname{\color[rgb]{1,0.3,0}}%
      \expandafter\def\csname LT8\endcsname{\color[rgb]{0.5,0.5,0.5}}%
    \else
      \def\colorrgb#1{\color{black}}%
      \def\colorgray#1{\color[gray]{#1}}%
      \expandafter\def\csname LTw\endcsname{\color{white}}%
      \expandafter\def\csname LTb\endcsname{\color{black}}%
      \expandafter\def\csname LTa\endcsname{\color{black}}%
      \expandafter\def\csname LT0\endcsname{\color{black}}%
      \expandafter\def\csname LT1\endcsname{\color{black}}%
      \expandafter\def\csname LT2\endcsname{\color{black}}%
      \expandafter\def\csname LT3\endcsname{\color{black}}%
      \expandafter\def\csname LT4\endcsname{\color{black}}%
      \expandafter\def\csname LT5\endcsname{\color{black}}%
      \expandafter\def\csname LT6\endcsname{\color{black}}%
      \expandafter\def\csname LT7\endcsname{\color{black}}%
      \expandafter\def\csname LT8\endcsname{\color{black}}%
    \fi
  \fi
    \setlength{\unitlength}{0.0500bp}%
    \ifx\gptboxheight\undefined%
      \newlength{\gptboxheight}%
      \newlength{\gptboxwidth}%
      \newsavebox{\gptboxtext}%
    \fi%
    \setlength{\fboxrule}{0.5pt}%
    \setlength{\fboxsep}{1pt}%
\begin{picture}(7936.00,3400.00)%
    \gplgaddtomacro\gplbacktext{%
      \csname LTb\endcsname
      \put(726,440){\makebox(0,0)[r]{\strut{}$0.7$}}%
      \put(726,897){\makebox(0,0)[r]{\strut{}$0.75$}}%
      \put(726,1353){\makebox(0,0)[r]{\strut{}$0.8$}}%
      \put(726,1810){\makebox(0,0)[r]{\strut{}$0.85$}}%
      \put(726,2266){\makebox(0,0)[r]{\strut{}$0.9$}}%
      \put(726,2723){\makebox(0,0)[r]{\strut{}$0.95$}}%
      \put(726,3179){\makebox(0,0)[r]{\strut{}$1$}}%
      \put(921,220){\makebox(0,0){\strut{}$0$}}%
      \put(1866,220){\makebox(0,0){\strut{}$250$}}%
      \put(2812,220){\makebox(0,0){\strut{}$500$}}%
      \put(3757,220){\makebox(0,0){\strut{}$750$}}%
      \put(4703,220){\makebox(0,0){\strut{}$1000$}}%
      \put(5648,220){\makebox(0,0){\strut{}$1250$}}%
      \put(6594,220){\makebox(0,0){\strut{}$1500$}}%
      \put(7539,220){\makebox(0,0){\strut{}$1750$}}%
    }%
    \gplgaddtomacro\gplfronttext{%
      \csname LTb\endcsname
      \put(6219,1053){\makebox(0,0)[l]{\strut{}quadratic}}%
      \csname LTb\endcsname
      \put(6219,833){\makebox(0,0)[l]{\strut{}degree 6}}%
      \csname LTb\endcsname
      \put(6219,613){\makebox(0,0)[l]{\strut{}degree 21}}%
    }%
    \gplbacktext
    \put(0,0){\includegraphics{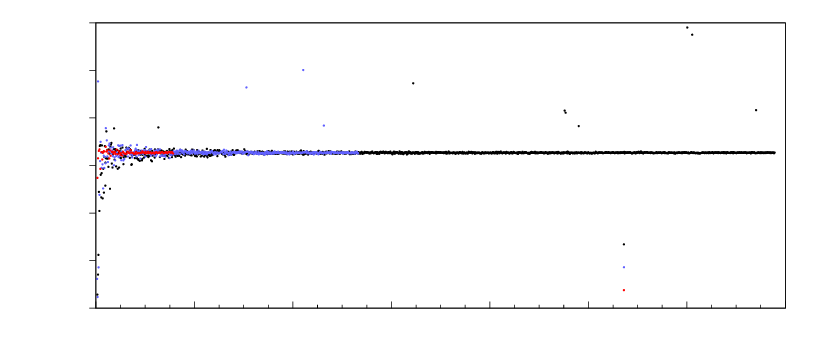}}%
    \gplfronttext
  \end{picture}%
\endgroup
\hspace{1cm}\mbox{.}\\
\caption{$\dfrac{T_{2+}(\K)}{T_2(\K)}$, in abscissa $\frac{1}{n}\log\Delta$.}
\label{fig:3A}
\end{fixedfig}
\smallskip
\noindent
Figure~\ref{fig:3A} shows that generally $T_{2+}(\K)$ is 13\% less than $T_2(\K)$, which is the current
result used by PARI/GP~\cite{PARI2}. However, according to the `no free-lunch' adage, the execution
time of the new algorithm approximatively doubles with respect to the previous one, but still remaining
around a few minutes for the worst case we have tested.
\smallskip

\begin{fixedfig}
\begingroup
  \makeatletter
  \providecommand\color[2][]{%
    \GenericError{(gnuplot) \space\space\space\@spaces}{%
      Package color not loaded in conjunction with
      terminal option `colourtext'%
    }{See the gnuplot documentation for explanation.%
    }{Either use 'blacktext' in gnuplot or load the package
      color.sty in LaTeX.}%
    \renewcommand\color[2][]{}%
  }%
  \providecommand\includegraphics[2][]{%
    \GenericError{(gnuplot) \space\space\space\@spaces}{%
      Package graphicx or graphics not loaded%
    }{See the gnuplot documentation for explanation.%
    }{The gnuplot epslatex terminal needs graphicx.sty or graphics.sty.}%
    \renewcommand\includegraphics[2][]{}%
  }%
  \providecommand\rotatebox[2]{#2}%
  \@ifundefined{ifGPcolor}{%
    \newif\ifGPcolor
    \GPcolorfalse
  }{}%
  \@ifundefined{ifGPblacktext}{%
    \newif\ifGPblacktext
    \GPblacktexttrue
  }{}%
  \let\gplgaddtomacro\g@addto@macro
  \gdef\gplbacktext{}%
  \gdef\gplfronttext{}%
  \makeatother
  \ifGPblacktext
    \def\colorrgb#1{}%
    \def\colorgray#1{}%
  \else
    \ifGPcolor
      \def\colorrgb#1{\color[rgb]{#1}}%
      \def\colorgray#1{\color[gray]{#1}}%
      \expandafter\def\csname LTw\endcsname{\color{white}}%
      \expandafter\def\csname LTb\endcsname{\color{black}}%
      \expandafter\def\csname LTa\endcsname{\color{black}}%
      \expandafter\def\csname LT0\endcsname{\color[rgb]{1,0,0}}%
      \expandafter\def\csname LT1\endcsname{\color[rgb]{0,1,0}}%
      \expandafter\def\csname LT2\endcsname{\color[rgb]{0,0,1}}%
      \expandafter\def\csname LT3\endcsname{\color[rgb]{1,0,1}}%
      \expandafter\def\csname LT4\endcsname{\color[rgb]{0,1,1}}%
      \expandafter\def\csname LT5\endcsname{\color[rgb]{1,1,0}}%
      \expandafter\def\csname LT6\endcsname{\color[rgb]{0,0,0}}%
      \expandafter\def\csname LT7\endcsname{\color[rgb]{1,0.3,0}}%
      \expandafter\def\csname LT8\endcsname{\color[rgb]{0.5,0.5,0.5}}%
    \else
      \def\colorrgb#1{\color{black}}%
      \def\colorgray#1{\color[gray]{#1}}%
      \expandafter\def\csname LTw\endcsname{\color{white}}%
      \expandafter\def\csname LTb\endcsname{\color{black}}%
      \expandafter\def\csname LTa\endcsname{\color{black}}%
      \expandafter\def\csname LT0\endcsname{\color{black}}%
      \expandafter\def\csname LT1\endcsname{\color{black}}%
      \expandafter\def\csname LT2\endcsname{\color{black}}%
      \expandafter\def\csname LT3\endcsname{\color{black}}%
      \expandafter\def\csname LT4\endcsname{\color{black}}%
      \expandafter\def\csname LT5\endcsname{\color{black}}%
      \expandafter\def\csname LT6\endcsname{\color{black}}%
      \expandafter\def\csname LT7\endcsname{\color{black}}%
      \expandafter\def\csname LT8\endcsname{\color{black}}%
    \fi
  \fi
    \setlength{\unitlength}{0.0500bp}%
    \ifx\gptboxheight\undefined%
      \newlength{\gptboxheight}%
      \newlength{\gptboxwidth}%
      \newsavebox{\gptboxtext}%
    \fi%
    \setlength{\fboxrule}{0.5pt}%
    \setlength{\fboxsep}{1pt}%
\begin{picture}(7936.00,3400.00)%
    \gplgaddtomacro\gplbacktext{%
      \csname LTb\endcsname%
      \put(594,440){\makebox(0,0)[r]{\strut{}$1.2$}}%
      \put(594,889){\makebox(0,0)[r]{\strut{}$1.3$}}%
      \put(594,1338){\makebox(0,0)[r]{\strut{}$1.4$}}%
      \put(594,1788){\makebox(0,0)[r]{\strut{}$1.5$}}%
      \put(594,2237){\makebox(0,0)[r]{\strut{}$1.6$}}%
      \put(594,2686){\makebox(0,0)[r]{\strut{}$1.7$}}%
      \put(594,3135){\makebox(0,0)[r]{\strut{}$1.8$}}%
      \put(789,220){\makebox(0,0){\strut{}$0$}}%
      \put(1753,220){\makebox(0,0){\strut{}$250$}}%
      \put(2718,220){\makebox(0,0){\strut{}$500$}}%
      \put(3682,220){\makebox(0,0){\strut{}$750$}}%
      \put(4646,220){\makebox(0,0){\strut{}$1000$}}%
      \put(5610,220){\makebox(0,0){\strut{}$1250$}}%
      \put(6575,220){\makebox(0,0){\strut{}$1500$}}%
      \put(7539,220){\makebox(0,0){\strut{}$1750$}}%
    }%
    \gplgaddtomacro\gplfronttext{%
      \colorrgb{0.58,0.00,0.83}%
      \put(6219,2962){\makebox(0,0)[l]{\strut{}quadratic}}%
      \colorrgb{0.00,0.62,0.45}%
      \put(6219,2742){\makebox(0,0)[l]{\strut{}degree 6}}%
      \colorrgb{0.34,0.71,0.91}%
      \put(6219,2522){\makebox(0,0)[l]{\strut{}degree 21}}%
    }%
    \gplbacktext
    \put(0,0){\includegraphics{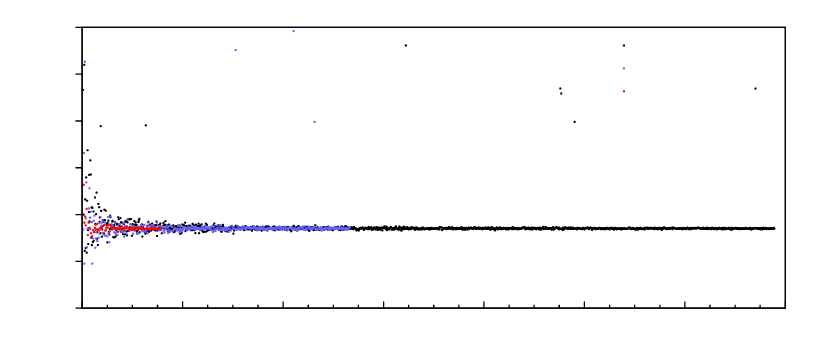}}%
    \gplfronttext
  \end{picture}%
\endgroup
\hspace{1cm}\mbox{.}\\
\caption{The resulting $c$, in abscissa $\frac{1}{n}\log\Delta$.}
\label{fig:4A}
\end{fixedfig}
\smallskip
\noindent
Figure~\ref{fig:4A} shows that the reduction in the size of the generating family comes from a value for
$c$ close to $1.37$ for all three series, in particular it is essentially independent of the degree of
the field but does not really depend on the discriminant either.
\smallskip

It is interesting to remark that both the findings in Figures~\ref{fig:2A} and~\ref{fig:4A} would follow
in case we would be able to improve two key points in our argument. Firstly, suppose that the field
satisfies a better version of Proposition~\ref{prop:1A}, namely
\[
\sum_{\substack{T< \Norm\p \leq cT}} \Lambda_{\K}(\p)\log(cT/\Norm\p)
\leq (c-1-\log c)T + O_n(\sqrt{T}\log^2T) + o(\sqrt{T}\log\Delta),
\]
hence without the factor $n$ as in Proposition~\ref{prop:1A}. In this case, the coefficient of
$\log\Delta$ for $\sqrt{T}$ becomes
\[
\frac{2\sqrt{c}}{c - 2(c-1-\log c)},
\]
whose minimum appears exactly for $c=1.37\ldots$, as we see in Figure~\ref{fig:4A}.
%
%
Moreover, the square of said minimum is $4\cdot 0.863\ldots$,
%
%
and since the factor $4$ comes from the need to bound $|\Ree\widehat{\psi^+}|$ by $|\widehat{\psi^+}|$
from~\eqref{eq:7A} to~\eqref{eq:8A}, in case we would be able to improve this step we would get exactly
the constant $0.863\ldots$ that we observe numerically in Figure~\ref{fig:2A}. %


\end{document}